\documentclass{article}

\usepackage[english]{babel}
\usepackage{amsmath}
\usepackage{graphicx}
\usepackage[colorinlistoftodos]{todonotes}

\title{Number of unique Edge-magic total labelings on Path $P_n$ }

\author{ Mukkai S. Krishnamoorthy\\ Rensselaer Polytechnic Institute, Troy, NY \\ Allen Lavoie \\ Washington University, St. Louis, MO\\ Ali Dasdan,\\ Hypergrowth, San Fransico Bay Area, CA\\ \and Bharath Santosh \\ Rennselaer Polytechnic Institute, Troy, NY  }

\date{\today}

\begin{document}
\maketitle

\begin{abstract}
Edge-magic total labeling was introduced by \cite{GS}. The number of edge-magic solutions for cycles have been explored in \cite{BS}. This sequence is mentioned in On Line Encyclopedia of Integer Sequences (OEIS) \cite{oeis}. In this short note, we enumerate the number of unique edge-magic total labelings on Path $P_n$
\end{abstract}

\section{Introduction}

\begin{figure}
\centering
\includegraphics[width=0.5\textwidth]{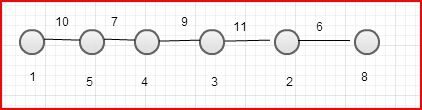}
\caption{\label{fig:fig1}{Path of length 5}}
\end{figure}

Edge-magic labeling (EMTL) has been studied in the past with an application towards communication networks. Given a simple undirected graph $G=(V,E)$, let $\lambda$ be a mapping from the numbers $1,2,\cdots, |V|+|E|$ to the vertices and edges of graph G, such that each element has an unique label.
The weight of an edge is obtained as the sum of the labels of that edge and its two end vertices. An edge-magic total labeling is labeling in which the weight of every edge is the same. The weight of the each edge is said to be a magic constant. Figure~\ref{fig:fig1}  illustrates an example for a path of length 5 with a magic constant of 16.  

Paper by Wallis et al \cite {wallis} describes existence of edge-magic total labeling of many types of graphs including $P_n$. The aim of this note is to enumerate unique edge-magic total labeling of Path $P_n$.

\section{Main Results}
\label{mainresults}

Our results are summarized in the following Table.

\begin{table}
\centering
\begin{tabular}{||l|r||}
\hline\hline Path Length & Number of Solutions \\\hline
0 & 1 \\
1 & 3 \\
2 & 12 \\
3 & 28 \\
4 & 48 \\
5 & 240 \\
6 & 944 \\
7 & 5344 \\
8 & 23408 \\
9 & 133808 \\
10 & 751008 \\
11 & 5222768 \\
12 & 37898776 \\
13 &	292271304 \\\hline
\end{tabular}
\caption{\label{tab:counts}Number of Edge-Magic Total Labels.}
\end{table}

As far as we have seen this series 1,3,12,28,48,240,944,5344,23408 does not appear in  OEIS.

\subsection{Method and Program}

We started with a simple python program to obtain all edge-magic solutions of paths of lengths 2 to 17. Paths of length 2 means that there will be three vertices and 2 edges, a total of 5 graph elements.

\begin{verbatim}

import itertools
for j in range(5,27,2):
    x = range(1,j+1)
    sum2 = 0
    for a in itertools.permutations(x):
        x = list(a)
        sum1 = x[0]+x[1]+x[2]
        d = 1
        for i in range(2,j-2,2):
            if (sum1== x[i]+x[i+1]+x[i+2]):
                d = 1
            else:
                d = 0
                break
        if (d==1):
            if (a[0]<a[j-1]):
                #print x, sum1
                sum2 = sum2+1
    print j, "\t", sum2
    
\end{verbatim}

This program generated one permutation at a time, and checked for whether the assignment leads to an edge-magic labeling.

However it is too slow and we could not compute past the path length of 7. We further optimized our python code and utilized the bounds on magic sum, $k$, similar to the one used in the paper\cite{BS}.
Let $f(r)=\frac{r\times (r+1)}{2}$

$\frac{f(2\times n +1) +f(n-1)}{n}  ~\le ~k ~\le \frac{ 2 \times f(2 \times n + 1)  - f( n+2)}{n}$

With this improvement and a shortcircuit optimization, we are able to get up to a path length of 13. All our code is located in the following github location \url {https://github.com/allenlavoie/path-counting }.

We will like to point the total number of edge-magic solutions for paths form a strict (albeit weak) upper bound for total number of edge-magic solutions for cycles of the same length.

\end{document}